\documentclass[a4paper,11pt]{article}
\usepackage{cleveref}
\usepackage{a4wide}
\usepackage{amssymb,amsmath,latexsym,amsthm} \usepackage{xcolor}
\usepackage{enumerate} %
\usepackage{dsfont}
\usepackage{geometry}
\usepackage{graphicx}

\usepackage{authblk} 
\usepackage[showonlyrefs]{mathtools} 
\usepackage[normalem]{ulem} %
\usepackage[latin1]{inputenc} 
\usepackage{yfonts} %
\usepackage{autonum}

\newcommand{\psum}[1]{\sum_{x\in L}{\vphantom{\sum}}'}
  { \end{itemize} }
  { \end{enumerate} }

\usepackage{array}
\newcolumntype{L}{>{\displaystyle} l <{}}

\newcounter{pcounter} 
 
\def\ED{\end{document}}

\newcommand{\IDEE}[1]{}

\renewcommand{\em}[1]{{\EMMM{#1}}}

\newcommand{\BS}[0]{\backslash}

\newcommand{\HW}[1]{} 

\renewcommand{\AA}{{\mathcal A}}

\newcommand{\FF}{{\mathcal F}}

\newcommand{\R}{\mathbb{R}}
\newcommand{\N}{\mathbb{N}}
\newcommand{\Z}{\mathbb{Z}}

\newcommand{\alp}{\alpha}
\newcommand{\bet}{\beta}

\newcommand{\eps}{\epsilon}

\renewcommand{\phi}{\varphi}

\renewcommand{\iint}{\int\!\!\!\!\int}

\def\XXint#1#2#3{{\setbox0=\hbox{$#1{#2#3}{\int}$}
\vcenter{\hbox{$#2#3$}}\kern-.5\wd0}}

\usepackage{color}
\definecolor{verylightblue}{rgb}{0.95, 0.95, 0.95}  
\definecolor{lightblue}{rgb}{0.7, 0.7, 1}
\definecolor{dblue}{rgb}{0.1, 0.1, 0.9}

\definecolor{eqyellow}{rgb}{0.9375,0.8984,0.5469}
\definecolor{subeqyellow}{rgb}{1,0.9373,0.8353}

\definecolor{mygreen}{rgb}{0.3, 0.6, 0.3} 
\definecolor{verylightgreen}{rgb}{0.95, 0.95, 0.95} 
\definecolor{verydarkgreen}{rgb}{0, 0.5, 0}
\definecolor{darkgreen}{rgb}{0 0.5 0}
\definecolor{mydarkgreen}{rgb}{0, 0.5, 0} 

\definecolor{mybrown}{rgb}{0.85, 0.4, 0.3}
\definecolor{verylightbrown}{rgb}{0.98, 0.72, 0.58}
\definecolor{verydarkbrown}{rgb}{0.44, 0.26, 0.26}

\definecolor{orange}{rgb}{1, 0.5, 0.2}

\definecolor{BurntOrange}{rgb}{1,0.356,0}

\definecolor{mydarkred}{rgb}{1,0.086,0.255}

\definecolor{RoseVYDP}{rgb}{0.84,0.086,0.255}

\definecolor{dgreen}{rgb}{0, 0.8, 0.5}     

\definecolor{CanaryBRT}{rgb}{1,0.76,0.26}

\definecolor{cyan}{rgb}{0, 1, 1}
\definecolor{verylightgray}{rgb}{0.95, 0.95, 0.95}
\definecolor{verylightgray}{rgb}{0.95, 0.95, 0.95}
\definecolor{verylightred}{rgb}{1, 0.8, 0.78}
\definecolor{verylightyellow}{rgb}{0.99, 0.98, 0.5}

\newcommand{\ignore}[1]{{}}

\newlength{\bibitemsep}\newlength{\bibparskip}
\let\oldthebibliography\thebibliography \renewcommand\thebibliography[1]{
  \oldthebibliography{#1} \setlength{\parskip}{\bibitemsep}
  \setlength{\itemsep}{\bibparskip} }
\makeatletter \newcommand{\VEC}[2][r]{originaler Name Spvek \gdef\@VORNE{1}
  \left(\hskip-\arraycolsep \begin{array}{#1}\vekSp@lten{#2}\end{array}
    \hskip-\arraycolsep\right) } \def\vekSp@lten#1{\xvekSp@lten#1;vekL@stLine;}
\def\vekL@stLine{vekL@stLine} \def\xvekSp@lten#1;{\def\temp{#1}
  \ifx\temp\vekL@stLine \else \ifnum\@VORNE=1\gdef\@VORNE{0} \else\@arraycr\fi
  #1 \expandafter\xvekSp@lten \fi} \makeatother

\newtheorem{theorem}{Theorem}[section]
\newtheorem{definition}[theorem]{Definition}
\newtheorem{corollary}[theorem]{Corollary}
\newtheorem{proposition}[theorem]{Proposition}
\newtheorem{conjecture}[theorem]{Conjecture} 
\theoremstyle{definition} \newtheorem{remark}[theorem]{Remark}

\numberwithin{equation}{section} 
\def\XXint#1#2#3{{\setbox0=\hbox{$#1{#2#3}{\int}$}
    \vcenter{\hbox{$#2#3$}}\kern-.5\wd0}} 
 \allowdisplaybreaks

\renewcommand{\em}[1]{\underline{#1}}

\makeatletter \newcommand\avsuminner[2]{ {\sbox0{$\m@th#1\sum$}
    \vphantom{\usebox0} \ooalign{ \hidewidth \smash{\vrule
        height\dimexpr\ht0+1pt\relax depth\dimexpr\dp0+1pt\relax} \hidewidth\cr
      $\m@th#1\sum$\cr } } } 
      \makeatother
      \author[*]{\rm Laurent B\'{e}termin} \author[**]{\rm Hans Kn\"upfer}
\author[**]{\rm Florian Nolte} \affil[*]{\normalsize{Faculty of Mathematics, University of Vienna, Oskar-Morgenstern-Platz 1, 1090 Vienna, Austria \texttt{laurent.betermin@univie.ac.at}. ORCID id: 0000-0003-4070-3344}}
\affil[**]{Institute of Applied Mathematics and IWR, University of Heidelberg, Im Neuenheimer Feld 205, 69120 Heidelberg, Germany. \texttt{knuepfer@uni-heidelberg.de}, \texttt{f.nolte@uni-heidelberg.de} }

\begin{document}

\title{Note on crystallization for  alternating particle
  chains}
\date\today

\maketitle

\begin{abstract}
  We investigate one-dimensional periodic chains of alternate type of particles
  interacting through mirror symmetric potentials. The optimality of the
  equidistant configuration at fixed density -- also called crystallization
  -- is shown in various settings. In particular, we prove the crystallization
  at any scale for neutral and non-neutral systems with inverse power laws
  interactions, including the three-dimensional Coulomb potential. We also
    show the minimality of the equidistant configuration at high density for
    systems involving inverse power laws and repulsion at the origin.
  Furthermore, we derive a necessary condition for crystallization at high
  density based on the positivity of the Fourier transform of the interaction
  potentials sum.
\end{abstract}

\setcounter{tocdepth}{1}

\section{Introduction}

A fundamental question in the theory of crystallization is to understand why
many large systems of interacting particles exhibit the spontaneous formation of
periodic structures and how they can be explained by energy minimization
\cite{BlancLewin-2015}. Such periodic structures are observed in systems
consisting of identical particles but also appear in models composed of
different types of particles. For example, ionic compounds exhibit periodic
structures, even though different attractive and repulsive interaction
potentials between the ions are present \cite{Nanoworld}. In this paper, we
consider prototypical alternating chains of particles where particles of
different (resp. same) kind repel (resp. attract) each other at short distance
and investigate necessary and sufficient conditions for the optimality of
periodic (equidistant) configurations. This consideration is motivated for
instance by alternating chains of magnetic domain walls (see
e.g. \cite{KnuepferMuratovNolte-preprint}). We note that while
one-dimensional model systems do not occur commonly in nature, they can be
created by confinement (see e.g.  \cite{ionicchainnanotubes}).

\medskip

In this paper, once the charges are fixed, as well as the interaction between
species, we show the optimality of the equidistant configuration at fixed
density, among one-dimensional periodic configurations of alternating species in
different settings. The novelty of the paper consists in the systematic analysis
for the ground state energy of alternating two-particle systems. We assume
repulsive interaction at short distances between different species in order to
avoid a collapsing of the ground state. We will also show that in the neutral
Coulomb case or in the power-law case, the equidistant configuration is the
unique minimum of the energy, and we expect this result to hold for a large
class of interaction potentials leading to a new kind of \textit{universal
  optimality} as defined by Cohn and Kumar in \cite{CohnKumar} for
two-component systems (see Conjecture
\ref{conj:univopt}). We note that the model considered is chosen as a simple
prototype model. More general, it would be interesting to derive conditions for
periodicity in higher-dimensional systems, to consider systems of more than two
different particles or to consider the case of different ratio between the
involved species.

\medskip

For one-dimensional systems of identical particles, Ventevogel and Nijboer
\cite{VN1,VN3} have derived several results about the optimality of the
equidistant configuration. In their work, interacting potentials are mirror
symmetric and correspond to semi-empirical potentials used in molecular
simulations (see e.g. \cite[p. 624]{CondensMatter}). In particular, they proved
the optimality of the equidistant configuration for convex interaction
potentials and Lennard-Jones-type potentials (also called ``Mie potentials")
among periodic configurations.  A similar result by Radin \cite{Rad1} shows the
optimality of an equidistant configuration for the classical $(12,6)$
Lennard-Jones potential, when the number of points -- added alternatively to
both sides of the configuration -- goes to infinity, and a generalization of this result has been recently shown in \cite{JWST19} where the temperature is added. Another recent result by
Bandegi and Shirokoff \cite[Sect. 6.1]{Bandegi:2015aa} gives numerical evidences
for the global optimality of the equidistant configuration for some values of
the density and the parameters of the Morse potential using convex
relaxation.

\medskip

One-dimensional systems involving
power-laws and two kind of species have been numerically studied in
\cite{LeviMinarLespowerlaws} in a different perspective, changing the species
ratio and considering interaction only between species of the same kind. We note
that similar studies for binary mixtures or particles have been made in
dimension two for dipolar (inverse power law) interaction and the Yukawa
potential (see \cite{Assoudetalyukawa} and references therein).

\medskip

Systems with different types of particles also arise in other models such as
e.g. chains of interacting magnetic dipoles
(e.g. \cite{nanoselfassembly,Spinchainsorder}, see also Fig.
\ref{fig-dipoles}). Also the interaction of stripe type magnetic domains in thin
ferromagnetic films can be described in this setting, where the sign of the
interaction energy between two interfaces in this model depends on the number of
in-between interfaces so that this model can be viewed as a system of
alternating particles of two kinds (see
\cite{KnuepferMuratovNolte-preprint}). Furthermore, our work can be related to
classical models of spin chains, and more precisely to the works of Giuliani et
al. on spin lattice models for stripe formation \cite{GLL06,GS16}.  We
 note that the type of models investigated in this
paper might also be interesting for biological models related to swarming and
flocking between different species (although in a dynamical, higher dimensional
setting, see e.g.
\cite{ZoologyCarrillo,burgerattractrepuls,Bertozzitwospecies}).

\medskip

Furthermore, one-dimensional quantum models involving nuclei and electrons have
also been studied. Brascamp and Lieb \cite{BrascampLieb} as well as Aizenman and
Martin \cite{AizenmanMartin} gave important results on the optimality of the
equidistant configuration concerning the Jellium model, see also \cite{SandierSerfaty1d}. The emergence of crystallized state for
one-dimensional system embedded in a periodic energy landscape is furthermore
studied in \cite{StefFried19} by Friedrich and Stefanelli where it is shown
that equidistant ground states are generally not expected. Furthermore, Blanc
and Le Bris \cite{BlancLebris} proved the periodicity of the ground state for
the one-dimensional Thomas-Fermi-von-Weizs\"acker energy.

\medskip

For higher dimensional systems, partial progress has been made for special
potentials or in restricted settings, especially in two and three dimensions. In
the case of only one type of particle, most of the known results in two and
  three dimensions concern perturbations of hard-sphere potentials or special
oscillating potentials (see \cite{BlancLewin-2015,BDLPSquare} and references
therein). In the special 8- and 24-dimensional cases, crystallization has been
shown at fixed density for the class of completely monotone potentials in
\cite{CKMRV2Theta}. In the case of several type of interacting particles, the
existing results rely on the design of specific short-range potentials such that
the difference of repulsion at short distance forced the particles to be in a
certain configuration. This is the main difference between these works and the
present paper where the interaction potentials are long-range and standard in
molecular simulations. A first proof of crystallization was given by Radin in
\cite{Radinquasiper}. Using three different radially symmetric short-range
interaction potentials, he proved the minimality of a two-dimensional binary
quasiperiodic configuration.  Furthermore, Friedrich and Kreutz
\cite{Friedrich:2018aa,FrieKreutSquare} have shown the optimality of
two-dimensional structures with alternating charges. They explored the
possibilities to obtain a rock-salt structure or an alternate honeycomb lattice
as minimizer of an interaction energy involving only short-range radially
symmetric potentials.

\section{Setting and statement of main results}

\begin{figure}[!h]
  \centering
  \includegraphics[width=10cm]{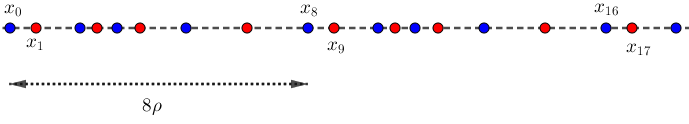}
  \caption{Example of periodic configuration $X\in \AA_{8}^\rho$.}
  \label{config}
\end{figure}

We consider a one-dimensional chain composed of two different types of
particles, located at the positions $x_k \in \R$ for $k \in \Z$ (see Fig.
\ref{config}). Since we are interested in the case when the interaction energy
between different types of particles is repulsive at short distance, we will assume that the
different particles are ordered alternatively. For technical
reasons, we also assume that the particles positions are $N$-periodic for a
  large even number $N$. We note that our results are independent of $N$ which is assumed to
  be very large. Recall that our goal is to show necessary conditions for
  minimizers of the considered systems to be $1$--periodic.
  \begin{definition}[Configurations and energy]  %
    Let $N \in 2\N$.
  \begin{enumerate}
  \item For $\rho>0$, the class of
    $N$--periodic configurations with density $\rho$ is denoted by
    \begin{align} \label{defxn} %
      \AA_{N}^\rho = \left\{ (x_k)_{k \in \Z} \ : \ \text{$x_0 = 0$,
      $x_k < x_{k+1}$ and $x_{k+N} = x_k +\rho^{-1} N$ $\forall k \in \Z$}
      \right\}.
    \end{align}
   We assume
      that the particle $x_k$ is of type $\eps_k$ where $\eps_k :=1$ if
      $k\in 2\Z+1$ and $\eps_k :=2$ if $k \in 2\Z$. The equidistant configuration
    $e^\rho \in \AA_N^\rho$ is denoted by
    $e^\rho := (k\rho^{-1})_{k \in \Z}$.
  \item For $\alp, \bet \in \{ 1, 2 \}$, let $f_{\alpha \beta} : \R \to \R$ be mirror
      symmetric interaction potentials, i.e
      $f_{\alp\bet}(-x)=f_{\alp\bet}(x)$. The associated energy for
      $\FF = (f_{11},f_{22},f_{12})$ is denoted by
  \begin{align}
    E_{\FF}(X):=&\frac{1}{N}\sum_{n=1}^N \sum_{k = -\infty\atop k\neq n}^\infty f_{\eps_n\eps_k}(x_{k}-x_n),%
                  \qquad \text{for $X = (x_i)_{i \in \Z} \in \AA_N^\rho$}.
  \end{align}
  \end{enumerate}
\end{definition}
We consider configurations that minimize the energy at fixed density. In
particular, we say that the equidistant configuration minimizes $E_\FF$ at any
scale if $e^\rho$ minimizes $E_\FF$ in $\mathcal{A}_N^\rho$ for all $\rho>0$ and
all $N\in 2\N$. Furthermore, we say that the equidistant configuration minimizes
$E_\FF$ at high density if there exists $\rho_0$ such that for all $N\in 2\N$
and all $\rho>\rho_0$, $e^\rho$ is a minimizer of $E_\FF$ in
$\mathcal{A}_N^\rho$.

\medskip

\begin{figure}[!h]
  \centering
  \includegraphics[width=10cm]{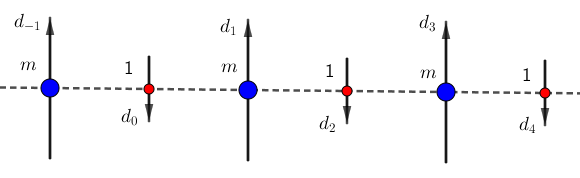}
  \caption{System of alternate oriented dipoles of magnitude $1$ and
    $m$.}
  \label{fig-dipoles}
\end{figure}

In \cite{VN1}, Ventevogel proved the optimality of $e^\rho$ for identical
  particles when, for any $\alpha,\beta\in \{1,2\}$, $f_{\alpha \beta}=f$ is a convex function. The
following theorem generalizes this result for two kinds of alternating species
and three kinds of interactions, also including some classes of nonconvex
functions:
\begin{theorem}[Sufficient condition]\label{thm-2}
  Suppose that $f_{\alpha \beta}(x)=\Phi_{\alpha \beta}^+(|x|)-\Phi_{\alpha \beta}^-(|x|)$ for convex and
  strongly tempered functions $\Phi_{\alpha \beta}^\pm : [0,\infty) \to \R$ for any
  $\alpha, \beta\in \{1,2\}$ in the sense that there exists $r_0,C,\eta>0$ such that
  \begin{align} 
    |\Phi_{\alpha \beta}^\pm (r)| < C r^{-1-\eta} \qquad \text{ for $r >r_0$.}
  \end{align}
  Moreover, suppose that the function $F$ is convex on $(0,+\infty)$, where
  \begin{align}\label{def-F}
    F(r) :=2\Phi_{12}^+(r)- \sum_{k=1}^\infty \big( \Phi_{12}^-((2k-1)r) + \Phi_{22}^-(2kr) + \Phi_{11}^-(2kr) \big).
  \end{align}
  Then the equidistant configuration is the unique minimizer of $E_\FF$
  at any scale.
\end{theorem}
We now give a direct application of Theorem \ref{thm-2} for systems with
alternating charges $1,-m$ and power-law interaction potential, in the
integrable case. We can think about two kinds of individuals with ``mass'' $1$
and $m$ interacting via the potentials $x\mapsto \pm |x|^{-p}$. The
following result shows that once $p$ is fixed, there exists an interval of
$m$ containing $m=1$ such that the equidistributed configuration is
the only minimizer of the energy at any scale. It gives (non-optimal) bounds on
$m$ such that this equilibrium is achieved.  Another physical motivation
are chains of antiparallel dipoles which are common e.g. in the self-assembly of
magnetic nanoparticles \cite[Fig. 1]{nanoselfassembly} and classical models of
spin chains \cite[Sect. 3]{Spinchainsorder} where this regular structure reaches
the equilibrium when there is no anisotropy field. Let us consider the following
toy model of a chain of dipoles $d_{n}$, located at position
$(x_{n},0,0) \in \R^3$ for $n \in \Z$. The dipoles are aligned in direction of
the $x_2$ axis with alternating orientation and with magnitude given by
$|d_{2k}| = 1$ and $|d_{2k+1}|=m$ (see Fig. \ref{fig-dipoles}). The
interaction potentials, up to a positive constant, are then given by
$f_{12}(x)=-m |x|^{-3}$, $f_{11}(x)=|x|^{-3}$, $f_{22}(x)= m^2 |x|^{-3}$. Hence, the following result gives a condition on
$m$ and $p$ such that the equidistant configuration is the only maximum for
this system at any scale.

\begin{corollary}[Crystallization for the inverse power-law: non-neutral
  case]\label{cor-integrablepowerlaw} \text{} %
  Let $f_{12}(x)=m |x|^{-p}$, $f_{11}(x)=-|x|^{-p}$ and
  $f_{22}(x)=-m^2 |x|^{-p}$, for $p>p_1$ where $p_1\approx 1.46498>1$ is the
  unique solution of $\zeta(p_1)=2^{p_1}$, where
    $\zeta(s):=\sum_{n>0} n^{-s}$ is Riemann's zeta function, and let $m$
  be such that
  \begin{align}\label{ineq-alpha}
    m_p< m < \frac{1}{m_p}, \quad m_p:=\frac{2^p -\sqrt{4^p - \zeta(p)^2}}{\zeta(p)}.
  \end{align}
  Then the equidistant configuration is the unique minimizer of $E_\FF$ at any scale.
\end{corollary}
We also note that if $m$ is sufficiently large (depending on any fixed $p$,
$N$) then the equidistributed configuration cannot be a minimizer of $E_\FF$,
the dominant interaction being $f_{22}(x)=-m^2 |x|^{-p}$, which forces the
particles to be close to each other.

\medskip

As a consequence of Theorem \ref{thm-2}, the Riesz potentials
$f_{11}(x)= f_{22}(x) = - f_{12}(x) = |x|^{-p}$ are minimized by the equidistant
configuration for any $p>1$. The next result improves this result to the case
$p \ge p_0$ where $p_0 \approx 0.655$ is the unique solution of
\begin{align} \label{def-prho} %
    \zeta(1+p_0)+1=2^{1+p_0} \qquad \text{ in $(0,\infty)$}.
\end{align}
In particular, this includes the case of the three-dimensional Coulomb potential.
\begin{theorem}[Riesz potentials]\label{thm-3}
  Let
  \begin{align}
    f_{12}(x)= - f_{11}(x)= - f_{22}(x) = \frac{1}{|x|^p} \qquad \text{for $p\ge p_0$, }
  \end{align}
  where $p_0$ is the unique solution of \eqref{def-prho}.  Then the equidistant
  configuration is the unique minimizer of $E_\FF$ at any scale.
\end{theorem}
\begin{remark}[Non-optimality at any scale for different exponents]
If we consider the anisotropic system where 
$$
f_{12}(x)=\frac{1}{|x|^p},\quad \textnormal{and}\quad f_{11}(x)=f_{22}(x)=-\frac{1}{|x|^q},\quad p\neq q,
$$
the crystallization on $e^\rho$ does not hold at any scale, contrary to what happens if $p=q$ as stated in Theorem \ref{thm-3}. Indeed, if $p>q$ (resp. $p<q$), the fact that $|x|^{-p}=o(|x|^{-q})$ as $|x|\to \infty$ (resp. as $|x|\to 0$) implies that there exists a density $\rho_0$ (resp. $\rho_1$) such that for all
$\rho<\rho_0$ (resp. $\rho>\rho_1$), $e^\rho$ is not a minimizer of $E_\FF$ in $\mathcal{A}_N^\rho$ since the main term of the energy is then attractive for this range of densities.
\end{remark}
We notice that the alternation of species is the only case where the minimizer
does not have two points  at the same location. Indeed, if
two points of the same kind are adjacent, it is sufficient to merge them in
order to get an energy equal to $-\infty$. In order to improve Theorem
\ref{thm-2} for non-summable potentials, the use of the homogeneity is a key
point.  In particular, Theorem \ref{thm-3} shows the maximality of the alternate
equidistant configuration (with charges $\pm 1$), at any scale, for the standard
Coulomb energy. Notice also that Theorem \ref{thm-3} is reminiscent of a result by Hubbard \cite{Hubbard} and Pokrovsky-Uimin \cite{PU78} on lattice spin systems (see also \cite{PB82}) but where the setting is slightly different: electrons are interacting through a Coulomb (or Riesz) interaction but a certain number of them are prescribed in a period and their presence at lattice sites are discussed. In a sense, there is no competition of two energies with totally opposite behavior in the latter as it is the case in our model.

\begin{remark}[\textbf{Limit of our method and universal optimality}]\label{rmk:univopt}
 We believe
  (supported by numerical experiments) that the equidistant configuration should
  be optimal for any $p > 0$ (including even the logarithmic
  potential). However, the constraint on the parameter $p_0$, is related to out
  method. We also believe that a new type of universal optimality -- in the
  sense of Cohn and Kumar \cite{CohnKumar} -- holds for such alternate systems,
  which leads to the following conjecture.
\begin{conjecture}[Universal optimality for alternating chains]\label{conj:univopt}
For any $f$ such that $f(r)=F(r^2)$ and $F$ is a completely monotone function (i.e. the Laplace transform of a nonnegative Borel measure), $e^\rho$ is the unique minimizer of $E_\FF$, where $\FF:=(-f,-f,f)$, at any scale.
\end{conjecture}  
Theorem \ref{thm-3} already proves this fact for a subclass of completely
monotone potentials. However, using the linear programming bounds method as in
\cite{CohnKumar,CKMRV2Theta} seems to be challenging -- even in dimension 1 --
for this type of alternate chain of particles. 
Furthermore, we expect Conjecture \ref{conj:univopt} to hold in higher dimensions $d\in \{2,4,8,24\}$ for the respective best packing and a more general class
of functions.
\end{remark}

\medskip

In order to explore a physically relevant anisotropic case, we now assume that the particles of the same type
  interact through the purely repulsive potential $f_{11}=f_{22}$ while the
  opposite type of particles are interacting through the attractive-repulsive
  one-well potential $f_{12}$, respectively given by
\begin{equation}\label{eq:potentials}
  f_{11}(x)=f_{22}(x) =\frac{1}{|x|^p}+\frac{1}{|x|^q}, %
  \quad f_{12}(x) = \frac{1}{|x|^p}-\frac{1}{|x|^q},\quad p>q>1.
\end{equation}
This model is
known as a good model for ionic interactions, because it takes into consideration the
Pauli exclusion principle by being repulsive at the origin, which avoids the
collapsing of the ground state. We therefore show, combining Theorem \ref{thm-2}
and \ref{thm-3}, that the equidistant configuration is the unique minimizer of
$E_{\mathcal{F}}$ among configurations of sufficiently high density.
\begin{theorem}[Optimality at high density in an anisotropic case]\label{thm_4}
Let $f_{ij}$ be defined by \eqref{eq:potentials}. Then there exists $\rho_0$ such that for all $\rho>\rho_0$ and all $N\in \N$, $e^{\rho}$ is the unique minimizer of $E_\FF$ in $\AA_N^\rho$. Furthermore, if $\rho<\rho_0$, then $e^\rho$ is not a minimizer of $E_\FF$ in $\AA_N^\rho$.
\end{theorem}
Notice that we have derived an analogous result in \cite{BFK20} in
any dimension in the case of (orthorhombic) Bravais lattices, showing the
minimality of the rock-salt structure. Furthermore, this type of high density result is usually a consequence of the universal optimality of a structure, as it is done for instance in \cite{BetTheta15,BFK20} and we think that Conjecture \ref{conj:univopt} could lead to a general high density result including a Coulombian tail $q=1$ in the same spirit as the two-dimensional result \cite[Thm. 1.1]{BetTheta15}.

\medskip

Finally, we recall a necessary condition for the optimality of the equidistant
configuration at high density, which has been derived by Ventevogel and
Nijboer's result in a related setting of identical particles \cite{VN3}, based
on the notion of strongly tempered potentials (cf. also S\"uto \cite{Suto2}).
\begin{proposition}[Necessary condition for high density crystallization] \label{prp-1} %
  Suppose that the functions $f_{\alpha \beta} \in C^2(\R)\cap L^1(\R)$ are mirror
  symmetric and strongly tempered for any $\alpha,\beta\in \{1,2\}$, in the sense that there exists $r_0,C,\eta>0$ such that  
    \begin{align} \label{strongly-tempered} %
      |f_{\alpha \beta} (x)| < C|x|^{-1-\eta} \qquad \text{ for any $|x|>r_0$.}
    \end{align}
   If the equidistant configuration is a minimizer of $E_\FF$ at high density, then
  \begin{align}
    \widehat{f_{12}}(k)+\frac{1}{2}\big(\widehat{f_{11}}(k)+\widehat{f_{22}}(k)  \big)\geq 0 && 
                                                                                                \text{for all } k\in \R,
  \end{align}
  where $\widehat{f}(k)=\frac{1}{\sqrt{2\pi}}\int_\R f(x)\cos(kx)dx$ is the
  Fourier transform of $f$.
\end{proposition}
The positivity of the Fourier transform of interaction potential seems to play a
crucial role in optimal point configurations theory (see
e.g. \cite{CohnKumar,Suto2,CKMRV2Theta}). For instance, Nikos et al. have
shown in \cite{Nikos07} -- for a generalized exponential model with one type of
particles -- an equivalence between the emergence of clusters and the existence
of negative Fourier components of the interaction potential. Proposition
\ref{prp-1} suggests that results in the same direction could be obtained for
two-components systems.
\section{Proofs}

In the following proofs, we use the notation $\ell:=\rho^{-1}$.

\subsection{Proof of Theorem \ref{thm-2}}
In this proof, for convenience, we write $\Phi_{\alpha \beta}^\pm(x)$ instead of $\Phi_{\alpha \beta}^\pm(|x|)$. In view of the assumptions of the theorem, for any $X\in \AA_N^\rho$ we have,
\begin{align}
  E_\FF(X)=&\frac{2}{N}\sum_{n=1}^N \sum_{k=1}^\infty \Phi_{12}^+(x_{n+2k-1}-x_n) + \frac{2}{N}\sum_{j=1}^{N/2}\sum_{k=1}^\infty \Phi_{22}^+(x_{2j+2k}-x_{2j})\\
       &+\frac{2}{N}\sum_{j=1}^{N/2}\sum_{k=1}^\infty \Phi_{11}^+(x_{2j-1+2k}-x_{2j-1})-\frac{2}{N}\sum_{n=1}^N \sum_{k=1}^\infty \Phi_{12}^-(x_{n+2k-1}-x_n)\\
       &-\frac{2}{N}\sum_{j=1}^{N/2} \sum_{k=1}^\infty \Phi_{22}^-(x_{2j+2k}-x_{2j})-\frac{2}{N}\sum_{j=1}^{N/2} \sum_{k=1}^\infty \Phi_{11}^-(x_{2j-1+2k}-x_{2j-1} )\\
       &=: S_1+S_2+S_3-S_4-S_5-S_6.
\end{align}
We estimate the six expressions using convexity of the functions and
periodicity. By convexity of $\Phi_{12}^+$  and with the notation
$d_n := x_{n+1}-x_n$, we have by Jensen's inequality
\begin{align}
  S_1&\geq \frac{2}{N}\sum_{n=1}^N \Phi_{12}^+(d_n)+2\sum_{k=2}^\infty \Phi_{12}^+\big((2k-1)\ell\big).
\end{align}
By convexity of $\Phi_{22}^+$, we furthermore obtain, using Jensen's inequality,
\begin{align}
  S_2&\geq \sum_{k=1}^\infty \Phi_{22}^+\Big( \frac{2}{N}\sum_{j=1}^{N/2} x_{2j+2k}-x_{2j}\Big)=\sum_{k=1}^\infty\Phi_{22}^+\big(2k\ell\big), 
\end{align}
the same holds for $S_3$ replacing $\Phi_{22}^+$ by $\Phi_{11}^+$.  For the
terms $S_4$, $S_5$, $S_6$ with a negative sign, we decompose into
nearest-neighbours distances. By convexity of $\Phi_{12}^-$ we get, again by Jensen's inequality
\begin{align}
  S_4 &=\frac{2}{N}\sum_{n=1}^N \sum_{k=1}^\infty \Phi_{12}^- \Big( \sum_{m=1}^{2k-1} d_{n+m-1} \Big) 
  \leq \frac{ 2}{N}\sum_{n=1}^N \sum_{k=1}^\infty \Phi_{12}^-((2k-1) d_n).
\end{align}
Similarly, by convexity of $\Phi_{22}^-$ and $\Phi_{11}^-$, we obtain
\begin{align}
  S_5&=\frac{2}{N}\sum_{j=1}^{N/2} \sum_{k=1}^\infty \Phi_{22}^-\Big( \sum_{m=1}^{2k} d_{2j+m-1} \Big) %
       \leq \frac{2}{N}\sum_{j=1}^{N/2} \sum_{k=1}^\infty \frac{1}{2k}\sum_{m=1}^{2k}  \Phi_{22}^-(2k d_{2j+m-1}) \\
     &\leq \frac{1}{N} \sum_{n=1}^N \sum_{k=1}^\infty \Phi_{22}^-(2k d_n), \\
  S_6
     &=\frac{2}{N}\sum_{j=1}^{N/2} \sum_{k=1}^\infty \Phi_{11}^- \Big( \sum_{m=1}^{2k} d_{2j+m-2} \Big) 
  \leq\frac{1}{N} \sum_{n=1}^N \sum_{k=1}^\infty \Phi_{11}^-(2k d_n).
\end{align}
Combining all these inequalities and since $F$ given by \eqref{def-F} is convex, we
have, by Jensen's inequality,
\begin{align}\label{ineq-F}
  E_\FF(X) \ %
  &\geq \ %
  2\sum_{k=1}^\infty \Phi_{12}^+\left((2k+1)\ell\right)+\sum_{k=1}^\infty \Phi_{22}^+\left(2k\ell\right)+\sum_{k=1}^\infty \Phi_{11}^+\left(2k\ell\right) %
    +\frac{1}{N}\sum_{n=1}^N F(d_n) \\
  &\geq 2\sum_{k=1}^\infty \Phi_{12}^+\left((2k+1)\ell\right)+\sum_{k=1}^\infty \Phi_{22}^+\left(2k\ell\right)+\sum_{k=1}^\infty \Phi_{11}^+\left(2k \ell\right)+F\left(\ell\right) %
        =E_\FF(e^\rho), \notag
\end{align}
with equality if and only if $X=e^\rho$.

\subsection{Proof of Theorem \ref{thm-3}}

The main idea is to compare the interaction on distances $|x_i - x_j|$ where
$i-j$ is even with interactions of distances where $i-j$ is odd: We use the
convex combination
\begin{align}\label{eq:id-convex}
  x_{n+2k}-x_n = \frac{(2k-j)}{2k} \frac{2k(x_{n+2k}-x_{n+j})}{2k-j} + \frac{j}{2k}\frac{2k(x_{n+j}-x_{n})}{j},
\end{align}
which holds for all $1\le j \le k$. We set $f(x) := |x|^{-p}$ and use again
  the notation $d_n := x_{n+1}-x_n$. Inserting \eqref{eq:id-convex}
for $j=1$ into $f$ and exploiting convexity, we get
\begin{align}
  f(x_{n+2k}-x_n) \le \frac{2k-1}{2k}f\left(\frac{2k(x_{n+2k}-x_{n+1})}{2k-1}\right) + \frac{1}{2k}f(2kd_{n}).
\end{align}
Since $f$ is homogeneous of degree $-p$, the last line implies
\begin{align}\label{eq:f-ub}
  f(x_{n+2k}-x_n) \le \Big(\frac{2k-1}{2k}\Big)^{1+p}f\left(x_{n+2k}-x_{n+1}\right) + \Big(\frac{1}{2k}\Big)^{1+p}f(d_n).
\end{align}
Averaging \eqref{eq:f-ub} over $n$ and using periodicity of $X$, we get
\begin{align}
  \frac{1}{N} \sum_{n=1}^N f(x_{n+2k}-x_n) \le \frac{1}{N} \sum_{n=1}^N \Big(\Big(\frac{2k-1}{2k}\Big)^{1+p}f(x_{n+2k-1}-x_{n}) + \Big(\frac{1}{2k}\Big)^{1+p}f(d_n)\Big),
\end{align}
i.e. a bound on the interaction on even distances in terms of the interaction on
odd distances. Inserting this estimate into the energy $E_\FF$ yields the lower
bound
\begin{align}
  E_\FF(X)\geq \frac{2}{N}\sum_{n=1}^N\sum_{k=1}^{\infty} a_kf(x_{n+2k-1}-x_n),
\end{align}
where the coefficients $a_k$ are given by
\begin{align}
  a_k:=
  \begin{cases}
    1-\frac{1}{2^{1+p}}-\sum_{j=1}^\infty\left(\frac{1}{2j}\right)^{1+p} =1- 2^{-(1+p)}(\zeta(1+p)+1) & \text{for }k=1,\\
    1-\left(\frac{2k-1}{2k}\right)^{1+p} & \text{otherwise}.
  \end{cases}
\end{align}
Since $a_1:p\mapsto 1- 2^{-(1+p)}(\zeta(1+p)+1)$ is an increasing function on
$(0,\infty)$, $p_0$ is unique and $p\geq p_0$ implies that $a_k \ge 0$ for
all $k \ge 1$.  Applying Jensen's inequality and inserting
$\frac{1}{N}\sum_{n=1}^{N} (x_{n+2k-1}-x_n) = (2k-1)\ell$ yields the lower bound
\begin{align}
  E_\FF(X) \ge 2\sum_{k=1}^{\infty}a_kf\left((2k-1)\ell\right) = E_\FF(e^\rho),
\end{align}
which is strict unless $X=e^\rho$, corresponding to equality in Jensen's
inequality.

\subsection{Proof of Corollary \ref{cor-integrablepowerlaw}}

  We want to apply Theorem \ref{thm-2}. We have
  \begin{align}
    F(r) = \frac{2m}{r^{p}}-\frac{(m^2 +1)\zeta(p)}{2^p r^{p}}, &&
    F''(r) = \Big(-\frac{\zeta(p)}{2^p}m^2 +2m -\frac{\zeta(p)}{2^p} \Big)\frac{p(p+1)}{r^{p+2}}.
  \end{align}
  The discriminant of the polynomial
  $P_p(m):=-\frac{\zeta(p)}{2^p}m^2 +2m -\frac{\zeta(p)}{2^p}$ is
  $\Delta=4\big( 1-\frac{\zeta(p)^2}{2^{2p}} \big)$, is positive if and
  only if $\zeta(p)<2^p$, since $p\mapsto 2^p -\zeta(p)$ is increasing on
  $(1,+\infty)$. Then, $P_p(m)>0$ if and only if $m$ satisfies
  \eqref{ineq-alpha}, and the proof is completed.

  \subsection{Proof of Theorem \ref{thm_4}}

    We first remark that, for any $\rho>0$, any $N\in \N$ and any $X\in \mathcal{A}_N^\rho$, we have
\begin{align}
  E_\FF(X)=E_{\FF_p}(X)+E_{\FF_q}(X),\quad \FF_p:=\left(\frac{1}{|x|^p},\frac{1}{|x|^p}, \frac{1}{|x|^p}  \right),\quad \FF_q:=\left(\frac{1}{|x|^q},\frac{1}{|x|^q}, -\frac{1}{|x|^q}  \right).
\end{align}
By the homogeneity of $f_{ij}$ and for the rescaled configuration
$X^1=\rho^{-1}X\in \mathcal{A}_N^1$, we have
$E_\FF(X)=\rho^p E_{\FF_p}(X^1)+\rho^{q}E_{\FF_q}(X^1)$. By Theorem \ref{thm-3}
(resp. Theorem \ref{thm-2}), we know that $E_{\FF_q}(X^1)\leq E_{\FF_q}(e^1)$
(resp. $E_{\FF_p}(X^1)\geq E_{\FF_p}(e^1)$) with equality if and only if
$X^1=e^1$. Hence, $E_\FF(X)\geq E_\FF(e^\rho)$ for all $X\in \mathcal{A}_N^\rho$
if and only if
\begin{align} \label{rhsof} %
 \rho> \rho_0 \ := \ \sup_{X^1\in \mathcal{A}_N^1\atop X_1\neq e^1}\left(\frac{E_{\FF_q}(e^1)-E_{\FF_q}(X^1)}{E_{\FF_p}(X^1)-E_{\FF_p}(e^1)}\right)^{\frac{1}{p-q}}.
\end{align}
It remains to show that $\rho_0 < +\infty$. Let $\|\cdot\|$ be the Euclidean norm on one period, with total length $N$, of any configuration belonging to $\AA_N^1$. 
We notice that $X^1\mapsto E_{\FF_q}(e^1)-E_{\FF_q}(X^1)$ and $X^1\mapsto E_{\FF_p}(X^1)-E_{\FF_p}(e^1)$ are both continuous in $\mathcal{A}_N^1$ with respect to $\|\cdot\|$ and vanish if and only if $X^1=e^1$. We also notice that the fact that $|x|^{-q}=o(|x|^{-p})$ as $x\to 0$ implies that any coalescence of two points of $X^1$ gives a zero energy contribution since the denominator tends to infinity faster than the numerator. Therefore, the above quotient is continuous on $\AA_N^1\backslash \{e^1\}$.

\medskip

By strict minimality (resp. maximality) of $e^1$ for $E_{\FF_p}$
(resp. $E_{\FF_q}$) and the Taylor expansion of these energies in
$\mathcal{A}_N^1$ close to $e^1$ combined with the continuity of the energies
and the fact that this strict minimality still holds in the limit $N\to \infty$,
it follows that there exists a sufficiently small $\delta>0$ and a constant
$C_\delta>0$ such that the right hand side of \eqref{rhsof} is bounded by
$C_\delta$ for all $X^1$ with $\|X^1-e^1\|<\delta$.
\medskip

We now show that the above quotient is bounded above uniformly in $N$ when $\|X^1-e^1\|\geq \delta$, $X^1\in \AA_N^1$. Since we
have already studied the behaviour of the above quotient as $X^1$ is
sufficiently close to $e^1$, we only need consider the quotient
$-E_{\FF_q}/E_{\FF_p}>0$, since $E_{\FF_q}\leq E_{\FF_q}<0$ and $E_{\FF_p}>0$ on $\AA_N^1$, and to show that there exists $C>0$ such that
\begin{align}\label{eq:equivalentbound}
  E_{p,q}(X^1):=C E_{\FF_p}(X^1) + E_{\FF_q}(X^1)\geq 0, \qquad \quad \forall N, \quad\forall X^1 = \{x_i\}_{i\in \Z} \in \AA_N^1.
\end{align}
Using that $r^{-p}\geq r^{-q}$ if and only if $r\leq 1$, we get
\begin{align}
  E_{p,q}(X^1) \ %
  &= \ \frac{1}{2N}\sum_{i=1}^N \sum_{x_j \in X^1 \BS \{x_i\}} \left\{ \frac{C}{|x_i-x_j|^p} + \frac{(-1)^{i-j}}{|x_i-x_j|^q} \right\}\\
  &\geq \frac{1}{2N}\sum_{i=1}^N \sum_{x_j \in X^1\backslash \{x_i\}:i-j\in 2\Z+1 \atop |x_i-x_j| \leq 1} \frac{C-1}{|x_i-x_j|^p} %
    - \frac{1}{2N}\sum_{i=1}^N \sum_{x_j \in X^1:  i-j \in 2\Z+1 \atop |x_i-x_j| > 1} \frac{1}{|x_i-x_j|^q} \\
   &\geq \ \frac{1}{2N}\left( C-1 - C_2\right).
\end{align}
For the last estimate, we have used that the first sum contains at least one
  term by the pigeonhole principle and we have assumed that $C\geq 1$. For
  the estimate of the second term, we have used that the
double sum is uniformly bounded because
\begin{align}
  \sum_{i=1}^N \sum_{x_j \in X^1:  i-j \in 2\Z+1 \atop |x_i-x_j| > 1} \frac{1}{|x_i-x_j|^q} \ 
  &
    \leq \ C_1 \iint_{(x,y)\in \R^2: |x-y|>1} \frac{dxdy}{|x-y|^q}<\ C_2 \ < \ \infty,
\end{align}
since $q>1$, where the constants $C_1,C_2$ only depend on $p,q$. For $C := 1+C_2$ we get $E_{p,q}(X^1)\geq 0$, which completes
the proof.

\subsection{Proof of Proposition \ref{prp-1}}

We adapt the proof of \cite[Sect. 2]{VN3} to our case. By assumption, we have,
for any $N$ and any $0<\ell<\ell_0$, $E_\FF(X)\geq E_\FF(e^\rho)$ for any $X = (x_i)_{i \in \Z} \in \AA_N^\rho$.  We
choose in particular $x_n :=y_n+n\ell-\eps$ with
$y_n=\eps \cos\left(\frac{2\pi m n}{N} \right)$ for some small $\eps>0$ such
that $\eps<\ell/2$ and for some $m\in \Z$. With this choice, we have $x_0 = 0$
and $x_{i+N}-x_i=N\ell$ and $x_{i+1}-x_i>0$ for any $i\in \Z$, and hence
$X\in \AA_N^\rho$.  Using Taylor expansion, by minimality of the equidistant
configuration we hence get, for $\ell$ sufficiently small,
\begin{align}
  \sum_{n=1}^N\sum_{j \in \Z}^\infty \left| y_j-y_n \right|^2 f_{\eps_j\eps_n}''((j-n)\ell) \geq 0.
\end{align}
Hence there is $\ell_1>0$ such that for any
$\ell\leq \ell_1\leq \ell_0$ and $g(x):=\frac 12(f_{11}(x)+f_{22}(x))$,
\begin{align} \label{thisin} %
  \sum_{j \in \Z} \left( 1-\cos\left((2j-1)q \right) \right)
  f_{12}''((2j-1)\ell)+\sum_{j\in \Z}\left( 1-\cos\left( 2q j\right)
  \right) g''(2j\ell)\geq 0,
\end{align}
where $q := \frac {2\pi m}N$.  Since \eqref{thisin} holds independently of $N$,
by an approximation argument we also have for any $x\in \R$ and any
$0<\ell\leq \ell_1\leq \ell_0$,
\begin{align}\label{condx}
  \sum_{j \in \Z} (1-\cos(2j x-x)f_{12}''((2j-1)\ell)+\sum_{j\in \Z} (1-\cos (2j x)) g''(2j\ell)\geq 0.
\end{align}
Thus, multiplying \eqref{condx} by $\ell$, taking $x=\ell k$ and dividing by
$k^2$, we get
\begin{align}
  0 &\leq \lim_{\ell\to 0} \Big(\ell\sum_{j \in \Z} \frac{(1-\cos( (2j \ell -\ell)k)}{k^2}f_{12}''((2j-1)\ell)+\ell\sum_{j \in \Z} \frac{(1-\cos (2j k\ell))}{k^2}g''(2j\ell)\Big)\\
    &=\int_{ \R} \frac{1-\cos ((2y-1)k)}{k^2}f_{12}''(2y-1)dy+\int_{ \R} \frac{1-\cos (2yk)}{k^2}g''(2y)dy \\
    &= \widehat f_{12}(k) + \widehat g(k) \qquad \text{for all $k \in \R \BS \{ 0 \}$.  }
\end{align}

  \medskip
  
  \paragraph{Acknowledgement.} %
  LB is grateful for the support of MATCH during his stay in Heidelberg. He also
  acknowledges support from the ERC Advanced grant ``Mathematics of the
  Structure of Matter'' (no. 321029) and from VILLUM FONDEN via the QMATH Centre
  of Excellence (no. 10059) during his stay in Copenhagen, as well as the WWTF
  research project ``Variational Modeling of Carbon Nanostructures"
  (no. MA14-009) at University of Vienna. HK acknowledges
support from the Deutsche Forschungsgemeinschaft (DFG, German Research Foundation) under Germany's Excellence Strategy - EXC-2181/1 - 390900948. We are also grateful to the anonymous referees for their suggestions. \renewcommand{\em}[1]{\it{#1}}

  {\small \bibliographystyle{plain} %
    \bibliography{crystoned}}

\end{document}